\documentclass[12pt]{article} 
\usepackage{amssymb,latexsym,color,psfig}
\newcounter{environment}[section]
\renewcommand{\theenvironment}{%
\arabic{section}.\arabic{environment}}

{\begin{rm}\refstepcounter{environment}{\textbf\theenvironment\
\bf Definition.~~}}%
{\end{rm}}

{\begin{rm}\refstepcounter{environment}{\textbf\theenvironment\
\bf Example.~~}}%
{\end{rm}}

{\begin{rm}\refstepcounter{environment}{\textbf\theenvironment\
\bf Conjecture.~~}}%
{\end{rm}}

{\begin{rm}\refstepcounter{environment}{\textbf\theenvironment\
\bf Proposition.~~}}%
{\end{rm}}

{\begin{rm}\refstepcounter{environment}{\textbf\theenvironment\
\bf Proposition and Definition.~~}}%
{\end{rm}}

{\begin{rm}\refstepcounter{environment}{\textbf\theenvironment\
\bf Theorem.~~}}%
{\end{rm}}

{\begin{rm}\refstepcounter{environment}{\textbf\theenvironment\
\bf Corollary.~~}}%
{\end{rm}}

{\begin{rm}\refstepcounter{environment}{\textbf\theenvironment\
\bf Lemma.~~}}%
{\end{rm}}

\savebox0{\rule[-2\unitlength]{0pt}{10\unitlength}%
\begin{picture}(10,10)(0,2)
\put(0,0){\line(0,1){10}}
\put(0,10){\line(1,0){10}}
\put(10,0){\line(0,1){10}}
\put(0,0){\line(1,0){10}}
\end{picture}}
\savebox1{\rule[-2\unitlength]{0pt}{10\unitlength}%
\begin{picture}(10,10)(0,2)
\put(0,10){\line(1,0){10}}
\put(10,0){\line(0,1){10}}
\put(0,10){\line(1,-1){10}}
\end{picture}}

\newlength\cellsize 
\savebox2{%
\begin{picture}(18,18)
\put(0,0){\line(1,0){18}}
\put(0,0){\line(0,1){18}}
\put(18,0){\line(0,1){18}}
\put(0,18){\line(1,0){18}}
\end{picture}}
\newcommand\cellify[1]{\def\thearg{#1}\def\nothing{}%
\ifx\thearg\nothing
\vrule width0pt height\cellsize depth0pt\else
\hbox to 0pt{\usebox2\hss}\fi%
\vbox to 18\unitlength{
\vss
\hbox to 18\unitlength{\hss$#1$\hss}
\vss}}
\newcommand\tableau[1]{\vtop{\let\\=\cr
\setlength\baselineskip{-16000pt}
\setlength\lineskiplimit{16000pt}
\setlength\lineskip{0pt}
\halign{&\cellify{##}\cr#1\crcr}}}
\savebox3{%
\begin{picture}(15,15)
\put(0,0){\line(1,0){15}}
\put(0,0){\line(0,1){15}}
\put(15,0){\line(0,1){15}}
\put(0,15){\line(1,0){15}}
\end{picture}}
\newcommand\expath[1]{%
\hbox to 0pt{\usebox3\hss}%
\vbox to 15\unitlength{
\vss
\hbox to 15\unitlength{\hss$#1$\hss}
\vss}}
\begin{document}
\newcommand{\qbc}[2]{ {\left [{#1 \atop #2}\right ]}}
\newcommand{\anbc}[2]{{\left\langle {#1 \atop #2} \right\rangle}}
\newcommand{\bbc}[2]{{\left( \!\! \left( {#1\atop #2}
  \right)\!\!\right) }}
\newcommand{\be}{\begin{enumerate}}
\newcommand{\ee}{\end{enumerate}}
\newcommand{\beq}{\begin{equation}}
\newcommand{\eeq}{\end{equation}}
\newcommand{\bea}{\begin{eqnarray}}
\newcommand{\eea}{\end{eqnarray}}
\newcommand{\beas}{\begin{eqnarray*}}
\newcommand{\eeas}{\end{eqnarray*}}
\newcommand{\zz}{\mathbb{Z}}
\newcommand{\pp}{\mathbb{P}}
\newcommand{\nn}{\mathbb{N}}
\newcommand{\qq}{\mathbb{Q}}
\newcommand{\rr}{\mathbb{R}}
\newcommand{\bm}[1]{{\mbox{\boldmath $#1$}}}
\newcommand{\sn}{\mathfrak{S}_n}
\newcommand{\la}{\lambda}
\newcommand{\ds}{\displaystyle}
\newcommand{\tla}{\tilde{\lambda}}
\newcommand{\fs}{\mathfrak{S}}
\newcommand{\ptq}{p\times q}
\newcommand{\st}{\,:\,}
\newcommand{\tr}{\textcolor{red}}
\newcommand{\tb}{\textcolor{blue}}
\newcommand{\tg}{\textcolor{green}}
\newcommand{\tm}{\textcolor{magenta}}
\newcommand{\tbn}{\textcolor{brown}}
\newcommand{\tp}{\textcolor{purple}}
\newcommand{\tn}{\textcolor{nice}}
\newcommand{\tor}{\textcolor{orange}}

\definecolor{brown}{cmyk}{0,0,.35,.65}
\definecolor{purple}{rgb}{.5,0,.5}
\definecolor{nice}{cmyk}{0,.5,.5,0}
\definecolor{orange}{cmyk}{0,.35,.65,0}

\begin{centering}
\textcolor{red}{\Large\bf Irreducible Symmetric Group Characters}\\
\textcolor{red}{\Large\bf of Rectangular Shape}\\[.2in] 
\textcolor{blue}{Richard P. Stanley}\footnote{Partially supported by
  NSF grant \#DMS-9988459 and by the Isaac Newton Institute for
  Mathematical Sciences.}\\ 
Department of Mathematics\\
Massachusetts Institute of Technology\\
Cambridge, MA 02139\\
\emph{e-mail:} rstan@math.mit.edu\\[.2in]
\textcolor{magenta}{version of 17 December 2002}\\[.2in]

\end{centering}
\vskip 10pt
\section{The main result.}
The irreducible characters $\chi^\lambda$ of the symmetric group $\sn$
are indexed by 
partitions $\lambda$ of $n$ (denoted $\lambda\vdash n$ or
$|\lambda|=n$), as discussed e.g.\ in
\cite[{\S}1.7]{macd} or \cite[{\S}7.18]{ec2}. If $w\in \sn$ has cycle type
$\nu\vdash n$ then we write $\chi^\lambda(\nu)$ for
$\chi^\lambda(w)$. If $\lambda$ has exactly $p$
parts, all equal to $q$, then we say that $\lambda$ has
\emph{rectangular shape} and write $\lambda=\ptq$. In this paper
we give a new formula for the values of the character $\chi^{p\times
q}$. 

Let $\mu$ be a partition of $k\leq n$, and let $(\mu,1^{n-k})$ be the
partition obtained by adding $n-k$ 1's to $\mu$. Thus $(\mu,1^{n-k})
\vdash n$. Define the \emph{normalized character}
$\widehat{\chi}^\lambda(\mu,1^{n-k})$ by
  $$ \widehat{\chi}^\lambda(\mu,1^{n-k}) = \frac{(n)_k
    \chi^\lambda(\mu,1^{n-k})}{\chi^\lambda(1^n)}, $$
where $\chi^\lambda(1^n)$ denotes the dimension of the character
$\chi^\lambda$ and $(n)_k=n(n-1)\cdots (n-k+1)$. Thus
\cite[(7.6)(ii)]{macd}\cite[p.\ 349]{ec2} 
$\chi^\lambda(1^n)$ is the number $f^\lambda$ of standard Young
tableaux of shape $\lambda$. Identify $\lambda$ with its
\emph{diagram} $\{(i,j)\st 1\leq j\leq \lambda_i\}$, and regard the
points $(i,j)\in\lambda$ as squares (forming the Young diagram of
$\lambda$). We write diagrams in ``English notation,'' with the first
coordinate increasing from top to bottom and the second coordinate
from left to right.
Let $\lambda=(\lambda_1,\lambda_2,\dots)$ and $\lambda'=
(\lambda'_1,\lambda'_2,\dots)$, where $\lambda'$ is the conjugate
partition to $\lambda$. The \emph{hook length} of the square
$u=(i,j)\in\lambda$ is defined by 
   $$ h(u) = \lambda_i+\lambda'_j-i-j+1, $$
and the Frame-Robinson-Thrall \emph{hook length formula}
\cite[Exam.~I.5.2]{macd}\cite[Cor.~7.21.6]{ec2} states that
  $$ f^\lambda = \frac{n!}{\prod_{u\in\lambda}h(u)}. $$

For $w\in\sn$ let $\kappa(w)$ denote the number of cycles of $w$ (in
the disjoint cycle decomposition of $w$). The main result of this paper is
the following.

\textbf{Theorem 1.} \emph{Let $\mu\vdash k$ and fix a permutation
  $w_\mu\in \fs_k$ of cycle type $\mu$. Then}
$$ \widehat{\chi}^{\ptq}(\mu ,1^{pq-k})=(-1)^k\sum_{uv=w_\mu}
      p^{\kappa(u)}(-q)^{\kappa(v)}, $$ 
\emph{where the sum ranges over all $k!$ pairs
$(u,v)\in\fs_k\times\fs_k$ satisfying $uv=w_\mu$.}
 
The proof of Theorem~1 hinges on a combinatorial identity involving
hook lengths and contents. Recall \cite[Exam.\ I.1.3]{macd}\cite[p.\
373]{ec2} that the 
\emph{content} $c(u)$ of the square $u=(i,j)\in\lambda$ is defined by
$c(u)= j-i$. We write $s_\lambda(1^p)$ for the Schur function
$s_\lambda$ evaluated at $x_1=\cdots=x_p=1$, $x_i=0$ for $i>p$. A well
known identity \cite[Exam.\ I.3.4]{macd}\cite[Cor.\ 7.21.4]{ec2} in
the theory of symmetric functions asserts that
  \beq s_\lambda(1^p) = \prod_{u\in\lambda}\frac{p+c(u)}{h(u)}. 
   \label{eq:sop} \eeq
Since the right-hand side is a polynomial in $p$, it makes sense to
define
  \bea s_\lambda(1^{-q}) = \prod_{u\in\lambda}\frac{-q+c(u)}{h(u)}. 
     \label{eq:slmq} \eea
Equivalently,
$s_\lambda(1^{-q})=(-1)^{|\lambda|}s_{\lambda'}(1^q)$. Regard $p$ and 
$q$ as fixed, and let $\lambda=(\lambda_1,\dots,\lambda_p)\subseteq
\ptq$ (containment of diagrams). Define the partition
$\tilde{\lambda}=(\tilde{\lambda}_1,\dots,\tilde{\lambda}_p)$ by
  \beq \tilde{\lambda}_i = q-\lambda_{p+1-i}. \label{eq:tla} \eeq
Thus the diagram of $\tilde{\lambda}$ is obtained by removing from the
bottom-right corner of $\ptq$ the diagram of $\lambda$ rotated
$180^\circ$. Write 
  $$ H_\lambda = \prod_{u\in\lambda}h(u), $$
the product of the hook lengths of $\lambda$.

\textbf{Lemma.} \emph{With notation as above we have}
  $$ H_{\ptq} = (-1)^{|\la|}H_\la H_{\tla} s_\la(1^p)
              s_\la(1^{-q}). $$

\begin{figure}
 \centerline{\psfig{figure=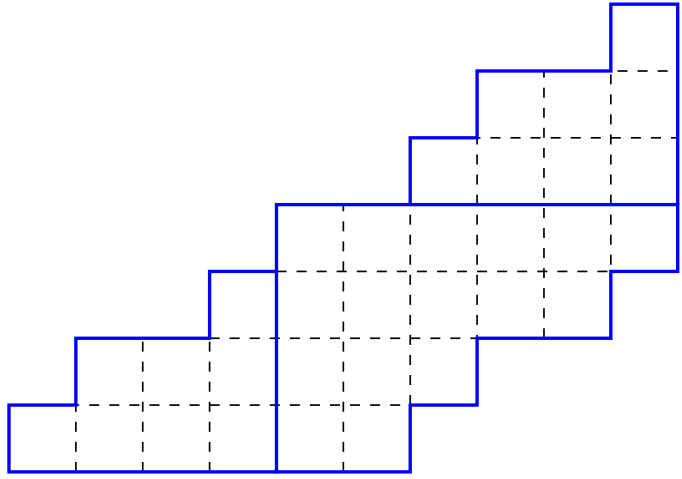}}
\caption{The shape SQ$(4,3,1)$ for $p=4$, $q=6$}
\label{fig:sq}
\end{figure}

\textbf{Proof.} Let $\lambda^\natural$ denote the shape $\lambda$
rotated $180^\circ$. Let SQ$(\lambda)$ denote the skew shape obtained
by removing $\lambda^\natural$ from the lower right-hand corner of
$\ptq$ and adjoining $\lambda^\natural$ at the right-hand end of the
top edge of $\ptq$ and at the bottom end of the left edge. See
Figure~\ref{fig:sq} for the case $p=4$, $q=6$, and $\lambda
=(4,3,1)$. It follows immediately from \cite[Thm.\ 1]{regev} that 
  \bea H_{\mathrm{SQ}(\lambda)} & = & H_{\tla}
     \prod_{u\in\lambda}(p+c(u)) \prod_{v\in\lambda'}(q+c(u))
   \nonumber\\ & = & (-1)^{|\lambda|}H_{\tla}
     \prod_{u\in\lambda}(p+c(u))(-q+c(u)). \label{eq:hsq1} \eea
It was proved in \cite{bes}\cite{janson}\cite{r-z} that the
multiset of hook lengths of the shape SQ$(\lambda)$ is the union of
those of the shapes $\ptq$ and $\lambda$, so
  \beq H_{\mathrm{SQ}(\lambda)} = H_{\ptq}H_\lambda. \label{eq:hsq2}
  \eeq
The proof now follows from equations (\ref{eq:sop}), (\ref{eq:slmq}),
(\ref{eq:hsq1}), and (\ref{eq:hsq2}). $\ \Box$
%

\textbf{Proof of Theorem~1.} Let $\ell=\ell(\mu)$. We first obtain an
expression for $\chi^{\ptq}(\mu ,1^{pq-k})$ using the
Murnaghan-Nakayama rule \cite[Exam.\ I.7.5]{macd}\cite[Thm.\
7.17.3]{ec2}. According to this rule,  
  $$ \chi^{\ptq}(\mu ,1^{pq-k}) = \sum_T (-1)^{\mathrm{ht}(T)}, $$
where $T$ ranges over all border-strip tableaux
$(B_1,B_2,\dots,B_{\ell+pq-k})$ of shape $\ptq$
and type $(\mu,1^{n-k})$. Here we are regarding $T$ as a sequence of
border strips removed successively from the shape $\ptq$. (See
\cite{macd} or \cite{ec2} for further details.) The first $\ell$
border strips $B_1,\dots,B_\ell$ will occupy some shape $\lambda\vdash
k$, rotated $180^\circ$, in the lower right-hand corner of $\ptq$. If
we fix this shape $\lambda$, then the number of choices for
$B_1,\dots,B_\ell$, weighted by $(-1)^{\mathrm{ht}(B_1)+\cdots+
\mathrm{ht}(B_\ell)}$, is by the Murnaghan-Nakayama rule just
$\chi^\lambda(\mu)$. The remaining border strips
$B_{\ell+1},\dots,B_{\ell+pq-k}$ all have one square (and hence height
0) and can be added in $f^{\tla}$ ways, where $\tla$ has the same
meaning as in (\ref{eq:tla}). Hence
  $$ \chi^{\ptq}(\mu,1^{pq-k}) = \sum_{{\lambda\subseteq \ptq
      \atop\lambda\vdash k}} \chi^\lambda(\mu)f^{\tla}, $$
so 
  \bea \widehat{\chi}(\mu,1^{pq-k}) & = & \frac{(pq)_k}{f^{\ptq}}
      \sum_{{\lambda\subseteq \ptq\atop\lambda\vdash k}}
     \chi^\lambda(\mu)f^{\tla} \nonumber\\ & = &
    \frac{(pq)_k H_{\ptq}}{(pq)!} \sum_{{\lambda\subseteq \ptq
    \atop\lambda\vdash k}}\chi^\lambda(\mu)\frac{(pq-k)!}{H_{\tla}}
     \nonumber\\ & = & H_{\ptq} \sum_{{\lambda\subseteq \ptq
      \atop\lambda\vdash k}}\chi^\lambda(\mu)H_{\tla}^{-1}.
       \label{eq:nchar} \eea
\indent Now let $\rho(w)$ denote the cycle type of a permutation
$w\in\fs_k$. The following identity appears in \cite[Prop.\
2.2]{h-s-s} and \cite[Exer.\ 7.70]{ec2}: 
  $$ \sum_{\lambda\vdash k}H_\lambda s_\lambda(x)s_\lambda(y)
     s_\lambda(z) = \frac{1}{k!}\sum_{{uvw=1\atop \mathrm{in}\ 
      \mathfrak{S}_k}}p_{\rho(u)}(x)p_{\rho(v)}(y)p_{\rho(w)}(z), $$
where $p_\nu(x)$ is a power sum symmetric function in the variables
$x=(x_1,x_2,\dots)$.
Set $x=1^p$, $y=1^{-q}$, take the scalar product (as defined in
\cite[{\S}I.4]{macd} or \cite[{\S}7.9]{ec2}) of both sides
with $p_\mu$, and multiply by $(-1)^k$. Since (in standard symmetric
function notation) the  
number of permutations in $\fs_k$ of cycle type $\mu$ is $k!/z_\mu$,
and since $\langle p_\mu,p_\mu\rangle=z_\mu$ and $\langle s_\lambda,
p_\mu\rangle =\chi^\lambda(\mu)$, we get
  \beq (-1)^k\sum_{\lambda\vdash k}H_\lambda s_\lambda(1^p)s_\lambda(1^{-q})
      \chi^\lambda(\mu) = (-1)^k\sum_{uv=w_\mu}p^{\kappa(u)}
       (-q)^{\kappa(v)}. \label{eq:sss} \eeq
Note that $s_\lambda(1^p)s_\lambda(1^{-q}) =0$ unless $\lambda
\subseteq \ptq$. Hence we can assume that $\lambda\subseteq \ptq$ in
the sum on the left-hand side of (\ref{eq:sss}). 

Now the coefficient of $\chi^\lambda(\mu)$ in (\ref{eq:nchar}) is
$H_{\ptq} H_{\tla}^{-1}$, while the coefficient of $\chi^\lambda(\mu)$
on the left-hand side of (\ref{eq:sss}) is $(-1)^k H_\lambda
s_\lambda(1^p)s_\lambda(1^{-q})$. By the lemma these two coefficients
are equal, and the proof follows. $\ \Box$

\section{Generalizations.} 
\begin{figure}
 \centerline{\psfig{figure=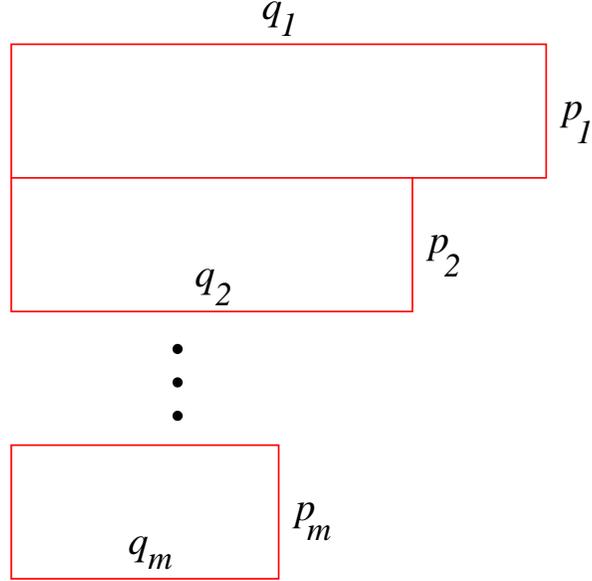}}
\caption{A union of $m$ rectangles}
\label{fig:nrec}
\end{figure}

The next step after rectangular shapes
would be shapes that are the union of two rectangles, then three
rectangles, etc. Figure~\ref{fig:nrec} shows a shape $\sigma\vdash
\sum_{i=1}^m p_iq_i$ that is a union of $m$ rectangles of sizes
$p_i\times q_i$, where $q_1>q_2>\cdots > q_m$. 

\textbf{Proposition 1.} \emph{Let $\sigma$ be the shape in
  Figure~\ref{fig:nrec}, and fix $k\geq 1$. Set $n=|\sigma|$ and}
  $$ F_k(p_1,\dots,p_m; q_1,\dots, q_m) =
  \widehat{\chi}^\sigma(k,1^{n-k}). $$
\emph{Then $F_k(p_1,\dots,p_m; q_1,\dots, q_m)$ is a polynomial
  function of the $p_i$'s and $q_i$'s with integer coefficients,
  satisfying} 
   $$ (-1)^k  F_k(1,\dots,1;-1,\dots,-1) = (k+m-1)_k. $$
\textbf{Proof.} Let $\lambda=(\lambda_1,\dots,\lambda_r)\vdash n$ and
  $$ \mu=(\mu_1,\dots,\mu_r)=(\lambda_1+r-1,\lambda_2+r-2,\dots,
     \lambda_r). $$
Define $\varphi(x)=\prod_{i=1}^r(x-\mu_i)$. A theorem of Frobenius
(see \cite[Exam.\ I.7.7]{macd}) asserts that
  \beq \widehat{\chi}^\lambda(k,1^{n-k}) = -\frac 1k [x^{-1}]
   \frac{(x)_k \varphi(x-k)}{\varphi(x)}, \label{eq:frob} \eeq
where $[x^{-1}]f(x)$ denotes the coefficient of $x^{-1}$ in the
expansion of $f(x)$ in \emph{descending} powers of $x$ (i.e., as a
Taylor series at $x=\infty$).

If we let $\lambda=\sigma$ in (\ref{eq:frob}) and cancel common factors
from the numerator and denominator, we obtain
  \bea \widehat{\chi}^\sigma(k,1^{n-k}) & = & -\frac 1k [x^{-1}]
  \frac{(x)_k\displaystyle \prod_{i=1}^m (x-(q_i+p_i+p_{i+1}+
   \cdots+p_m))_k}{\displaystyle \prod_{i=1}^m (x-(q_i+p_{i+1}
   +p_{i+2}+\cdots+p_m))_k} \label{eq:chires}\\ & = & 
    -\frac 1k[x^{-1}]H_k(x), \nonumber \eea
say. Since
  $$ \frac{1}{x-a} =\frac 1x+\frac{a}{x^2}+\frac{a^2}{x^3}+\cdots, $$
it is clear that $[x^{-1}]H_k(x)$ will be a polynomial
$F_k(p_1,\dots,p_m; q_1,\dots, q_m)$ in the $p_i$'s and $q_i$'s with
integer coefficients.  If we put $p_i=1$ and $q_i=-1$ then we obtain
(after cancelling common factors) 
  $$ F_k(1,\dots,1;-1,\dots,-1) = -\frac 1k [x^{-1}]
    \frac{(x-k+1)(x-m+1)_k}{x+1}. $$
Since the sum of the residues of a rational function $R(x)$ in the
extended complex plane is 0, it follows that 
 \beas -\frac 1k [x^{-1}]\frac{(x-k+1)(x-m+1)_k}{x+1} & = & -\frac 1k
   \mathrm{Res}_{x=-1}\left(\frac{(x-k+1)(x-m+1)_k}{x+1}\right)\\ & = &
   (-m)_k\\ & = & (-1)^k(k+m-1)_k. \eeas
It remains to show that the coefficients of  $F_k(p_1,\dots,p_m;
q_1,\dots, q_m)$ are integers. Equivalently, the coefficients of
the polynomial
  $$  [x^{-1}] \frac{(x)_k \varphi(x-k)}{\varphi(x)} $$
are divisible by $k$. But
  $$  \frac{(x)_k \varphi(x-k)}{\varphi(x)} \equiv (x)_k\
   (\mathrm{mod}\ k) $$
and 
  $$ [x^{-1}] (x)_k =0, $$
so the proof follows. $\ \Box$


\tm{\textsc{Note.}} For any fixed $\mu\vdash k$, J. Katriel has shown
(private communication), based on a method \cite{katriel} for
expressing $\widehat{\chi}^\lambda(\mu,1^{n-k})$ in terms of the values
$\widehat{\chi}^\lambda(j,1^{n-j})$, that $\widehat{\chi}^\sigma(\mu
,1^{n-k})$ is a polynomial $F_\mu(p_1,\dots,p_m;q_1,\dots,q_m)$ with
rational coefficients satisfying
  $$ (-1)^k  F_\mu(1,\dots,1;-1,\dots,-1) = (k+m-1)_k. $$
It can be deduced from the Murnaghan-Nakayama rule that in fact
the function $F_\mu(p_1,\dots,p_m;q_1,\dots,q_m)$ is a polynomial with
\emph{integer} coefficients. We conjecture that in fact the
coefficients of $F_\mu(p_1,\dots,p_m;q_1,\dots,q_m)$ are
nonnegative: 

\textbf{Conjecture 1.} \emph{For fixed $\mu\vdash k$,
  $\widehat{\chi}^\sigma(\mu ,1^{n-k})$ is a polynomial
  $F_\mu(p_1,\dots,p_m;$ $q_1,\dots,q_m)$ with integer coefficients such
  that $(-1)^kF_\mu(p_1,\dots,p_m;-q_1,\dots,-q_m)$ has nonnegative
  coefficients summing to $(k+m-1)_k$.}

We do not have a conjectured combinatorial interpretation of the
  coefficients of $(-1)^kF_\mu(p_1,\dots,p_m;-q_1,\dots,-q_m)$. When
  $m=2$ we have the following data, where we write $a=p_1$, $p=p_2$,
  $b=q_1$, $q=q_2$:
 \beas -F_1(a,p;-b,-q) & = & ab+pq\\
    F_2(a,p;-b,-q) & = & a^2b+ab^2+2apq+p^2q+pq^2\\
    -F_3(a,p;-b,-q) & = & a^3b+3a^2b^2+3a^2pq+ab^3+3abpq+3ap^2q+
             3apq^2\\ & & \ +p^3q+3p^2q^2+pq^3+ab+pq\\ \pagebreak
    F_4(a,p;-b-q) & = & a^4b+6a^3b^2+4a^3pq+6a^2b^3+12a^2bpq+6a^2p^2q\\ 
     & & \ + 6a^2pq^2+ab^4+4ab^2pq+4abp^2q+4abpq^2+4ap^3q\\ & & \
        +14ap^2q^2+4apq^3+p^4q+6p^3q^2+6p^2q^3+pq^4+5a^2b\\
        & & \ +5ab^2+10apq+5p^2q+5pq^2. \eeas
We can say something more specific about the leading terms of
$F_k(p_1,\dots,p_m;$ $q_1,\dots,q_m)$. Let $G_k(p_1,\dots,p_m;q_1,
  \dots,q_m)$ denote these leading terms, viz., the terms of total
  degree $k+1$. 

\textbf{Proposition 2.} \emph{We have}
  $$ \frac 1x+\sum_{k\geq 0} G_k(p_1,\dots,p_m;q_1,\dots,q_m)x^k = $$
    \beq \quad\frac{1}{\left( \frac{\ds x\prod_{i=1}^m
          (1-(q_i+p_{i+1}+p_{i+2}+\cdots 
         +p_m)x)}{\ds\prod_{i=1}^m (1-(q_i+p_i+p_{i+1}+\cdots+
          p_m)x)}\right)^{\langle -1\rangle}}, \label{eq:lt} \eeq
\emph{where $^{\langle -1\rangle}$ denotes compositional inverse
\cite[{\S}5.4]{ec2} with respect to $x$. In particular, the generating
function $\sum G_kx^k$ is algebraic over
$\mathbb{Q}(p_1,\dots,p_m,q_1,\dots,q_m,x)$.}  

\textbf{Proof.} From (\ref{eq:chires}) we have
  \beas G_k(p_1,\dots,p_k;q_1,\dots,q_k) & = & -\frac 1k [x^{-1}]
  \frac{x^k\ds\prod_{i=1}^m (x-(q_i+p_i+p_{i+1}+\cdots+p_m))^k}
       {\ds\displaystyle\prod_{i=1}^m
      (x-(q_i+p_{i+1}+p_{i+2}+\cdots+p_m))^k}\\ & = & 
          -\frac 1k[x^{-1}]L(x)^k, \eeas
say. Let $L(1/x)=M(x)/x$, so $M(0)=1$. Regard $M(x)$ as a power series
in ascending powers of $x$, i.e., an ordinary Taylor series at
$x=0$. Then by the Lagrange inversion
formula \cite[Thm.\ 5.4.2]{ec2} we have
  $$ [x^{-1}]L(x)^k =[x^{k+1}]M(x)^k=-k[x^k]\frac{1}{(x/M(x))^{\langle
      -1\rangle}}, $$ 
so equation (\ref{eq:lt}) follows. $\ \Box$

Proposition~2 was also proved by Philippe Biane (private
communication) in the same way as here, though using the language of
free probability theory.

It follows from Proposition~1 or Proposition~2 that $(-1)^k
G_k(p_1,\dots,p_m;$ $-q_1,\dots,-q_m)$ is a polynomial with integer
coefficients summing to 
  $$ S_k:=(-1)^kG_k(1,\dots,1;-1,\dots,-1). $$ 
From Proposition~2 we have
  $$ -\frac 1x+\sum_{k\geq 0}S_kx^k = \frac{-1}{\ds\left(
       \frac{x(1-x)}{1-(m-1)x}\right)^{\langle -1\rangle}}, $$
an algebraic function of degree two. When $m=1$ we have $S_k=C_k$, the
$k$th Catalan number. Hence by Theorem~1 $C_k$ is equal to the number
of pairs $(u,v)\in\fs_k\times\fs_k$ such that $\kappa(u)+\kappa(v) =
k+1$ and $uv=(1,2,\dots,k)$, a known result (e.g., \cite[Exer.\
6.19(hh)]{ec2}). Moreover, it follows easily from Proposition~2 that
  $$ (-1)^kG_k(p;-q) = \sum_{i=1}^k N(k,i)p^{k+1-i}q^i, $$
where $N(k,i)=\frac 1k{k\choose i}{k\choose i-1}$, a \emph{Narayana
  number} \cite[Exer 6.36]{ec2}. Hence $N(k,i)$ is equal to the number
of pairs $(u,v)\in\fs_k\times\fs_k$ such that $\kappa(u)=i$,
$\kappa(v)=k+1-i$, and $uv=(1,2,\dots,k)$.
When $m=2$ we have $S_k=r_k$, a (big) \emph{Schr\"oder number}
\cite[p.\ 178]{ec2}. 


It would follow from Conjecture 1 that the polynomial $(-1)^k
G_k(p_1,\dots, p_m;$ $-q_1,\dots,-q_m)$ has nonnegative coefficients. In
fact, Sergi Elizalde has shown (private communication of May, 2002)
that 
$$ \hspace{-2in}(-1)^k G_k(p_1,\ldots,p_m;-q_1,\ldots,-q_m) $$ 
$$ \quad ={1\over
k}\sum_{i_1+\cdots+i_m+j_1+\cdots+j_m=k+1}{k\choose i_1}\bbc{i_1}
{j_1} $$
 $$ \qquad \prod_{s=2}^m \left(\sum_{r=0}^{\min(i_s,j_s)}{k\choose
         r}\bbc{r} 
{j_s-r}{k-r-i_1-\cdots-i_{s-1}-j_1-\cdots-j_{s-1}\choose
i_s-r_s}\right) $$ 
  $$ \qquad\qquad\qquad {p_1^{i_1}\cdots p_m^{i_m}q_1^{j_1}\cdots
q_m^{j_m}},$$
where $\bbc{a}{b}={a+b-1\choose b}$.
Thus in particular $(-1)^k G_k(p_1,\dots, p_m; -q_1,\dots,-q_m)$
indeed does have nonnegative coefficients. Do they have a simple
combinatorial interpretation?


\end{document}